\documentclass[10pt,oneside]{jm}
\usepackage{amsfonts,amssymb,amsmath,rotating}
\usepackage[backref, pageanchor=false, colorlinks=true,%
            linkcolor=blue, citecolor=red, urlcolor=magenta]{hyperref}

\newcommand{\abs}{\vskip 0.5em\noindent\rm}
\newcommand{\Abs}{\paragraph{}\hspace{-1em}\rm}
\newcommand{\AbsT}[1]{\paragraph{\hspace{-1em} #1}\rm}

\newcommand{\bfi}{\noindent {\bf i)} }

\newcommand{\bfii}{\noindent {\bf ii)} }

\newcommand{\bfiii}{\noindent {\bf iii)} }

\newcommand{\bfa}{\noindent {\bf a)} }
\newcommand{\bfb}{\noindent {\bf b)} }

\newcommand{\absr}{\abs\hrulefill}

\newcommand{\F}{\mathbb F}

\newcommand{\N}{\mathbb N}

\newcommand{\Q}{\mathbb Q}

\newcommand{\Z}{\mathbb Z}

\newcommand{\cA}{\mathcal A}

\newcommand{\cE}{{\mathcal E}}

\newcommand{\ccH}{\mathcal H}

\newcommand{\cO}{\mathcal O}

\newcommand{\cQ}{\mathcal Q}

\newcommand{\cS}{\mathcal S}

\newcommand{\al}{\alpha}

\newcommand{\bt}{\beta}
\newcommand{\gm}{\gamma}

\newcommand{\dt}{\delta}

\newcommand{\om}{\omega}

\newcommand{\Om}{\Omega}

\newcommand{\ph}{\varphi}

\newcommand{\Aut}{\text{Aut}}

\newcommand{\End}{\text{End}}

\newcommand{\Fix}{\text{Fix}}

\newcommand{\Hom}{\text{Hom}}
\newcommand{\IBr}{\text{IBr}}

\newcommand{\Irr}{\text{Irr}}
\newcommand{\IPr}{\text{IPr}}

\newcommand{\rk}{\text{rk}}

\newcommand{\SL}{\text{SL}}

\newcommand{\Stab}{\text{Stab}}

\newcommand{\Tr}{\text{Tr}}

\newcommand{\ld}{,\ldots\hskip0em ,}
\newcommand{\lr}[1]{\langle #1\rangle}

\newcommand{\mt}{\mapsto}

\newcommand{\ra}{\rightarrow}

\newcommand{\Lra}{\leftrightarrow}

\newcommand{\cn}{\colon}

\newcommand{\sseq}{\subseteq}

\newcommand{\bsl}{\backslash}

\newcommand{\tm}{\times}

\newcommand{\GAP}{{\sf GAP}}

\newcommand{\MA}{{\sf MeatAxe}}

\newcommand{\ORB}{{\sf ORB}}

\begin{document}
\raggedbottom
\pagestyle{myheadings}
\markboth{}{}
\thispagestyle{empty}
\setcounter{page}{1}

\begin{center} \Large\bf
The projective cover of the trivial module \\ 
in characteristic $11$ for the sporadic \\
simple Janko group $J_4$ revisited
\vspace*{1em} \\
\large\rm Jürgen Müller \vspace*{1.5em} \\
\end{center}

\begin{abstract} \vspace*{-3em}\noindent\hrulefill\vspace*{2em} \\ \noindent 
This is a sequel to \cite{part1}, where we have determined the $11$-modular 
projective indecomposable summands of the permutation character of $J_4$ 
on the cosets of an $11'$-subgroup of maximal order, amongst them the 
projective cover of the trivial module, up to a certain parameter. Here,
we fix this parameter, by applying a new condensation method for induced
modules which uses enumeration techniques for long orbits.

\vspace*{0.5em} \noindent

\vspace*{0.5em} \noindent 

\vspace*{0em} \noindent \hrulefill 
\end{abstract}


\section{Introduction}\label{intro}

\Abs
The present article is a sequel to \cite{part1}, which is devoted to 
answering (to the negative) the question (posed in \cite{MalRob}) whether 
the projective cover of the trivial module in characteristic $11$ for 
the largest sporadic simple Janko group $G:=J_4$ is a permutation module.
Actually, this question boils down to determine whether or not the 
projective permutation character $1_H^G$, where $H<G$ is a maximal 
($11'$-)subgroup of shape $H\cong 2^{10}\cn L_5(2)$, is projective 
indecomposable. 

\abs
According to the {\sf ModularAtlasHomepage}, virtually nothing is known 
about the decomposition numbers of $G$ is characteristic $11$. Thus the
strategy in \cite{part1} was to find the decomposition of $1_H^G$ into 
projective indecomposable characters. It turned out that there are four 
distinct indecomposable summands, which are reproduced in Table \ref{pchtbl}
(on page \pageref{pchtbl}). This can be seen as the first step towards
the ambitious goal of determining the $11$-modular decomposition matrix 
of $G$, which is particularly compelling as $G$ has trivial-intersection 
Sylow $11$-subgroups.

\abs
Alone, in \cite{part1}, we have not been able to fix the parameter 
$a\in\{0,1\}$ appearing in Table \ref{pchtbl}. But this should be done
before proceeding to find more decomposition numbers of $G$. Hence the 
purpose of the present article is to close this gap, by showing that we 
actually have $a=0$, obeying to the conventional choices for decomposition 
maps made in the {\sf ModularAtlas}.

\abs
To this end, we invoke a maximal subgroup $U<G$ of shape
$U\cong U_3(11)\cn 2$, whose $11$-modular decomposition matrix is well-known. 
In order to relate the decomposition matrices of $U$ and $G$, it turns out 
that subtle details of the embedding of $U$ into $G$ play a crucial role 
here. These can be captured by a consideration of the automorphism group 
of the ordinary character table of $G$. 

\abs
Having this in place, letting $\cO$ be a $G$-set affording the permutation 
character $1_H^G$, and letting $\F$ be a field of characteristic $11$, 
we examine how the restriction of the permutation module $\F[\cO]$ to $U$ 
decomposes into projective indecomposable modules. To do so, we are finally
led to consider the action of $\End_{\F[G]}(\F[\cO])$ on
$\Hom_{\F[U]}(\F[\cO]|_U,V)\cong\Hom_{\F[G]}(\F[\cO],V^G)$,
for certain simple $\F[U]$-modules $V$.

\abs
The latter step essentially amounts to computing the `condensed module'
afforded by the module $V^G$ induced from the subgroup $U$, with respect 
to the `condensation subgroup' $H$.
(As a general reference for `condensation', see for example \cite{Habil}.)
A technique to compute condensed modules of induced modules, for subgroups 
$U$ of smallish index in $G$, has been developed in \cite{IndCond}.
The present approach, combining these ideas with the orbit enumeration
techniques available in the \GAP{} package \ORB{}, now allows both 
subgroups $U$ and $H$ to have large index in $G$.
We expect this to be of independent interest.

\Abs
The present article is organized as follows:
In the rest of Section \ref{intro} we indicate the computational 
tools we are using, and we sketch the ideas behind the orbit
enumeration techniques available in \ORB{}.
In Section \ref{indcond} we present the piece of theory underlying the
new condensation technique advertised above.
In Section \ref{chars} we collect some character theoretic facts on
$G$ and various of its subgroups, and we consider table automorphisms, 
to clarify where and where not choices can be made.
In Section \ref{enum} we enumerate $\cO$ by $U$-orbits, by applying \ORB{}.
In Section \ref{conclusion}, using an idea inspired by \cite{LMR}, 
we examine certain condensed induced modules, in order to finally 
determine the missing parameter.

\AbsT{Computational tools.}
To facilitate group theoretic and character theoretic computations we use
the computer algebra system \GAP{} \cite{GAP}, and its comprehensive database
{\sf CTblLib} \cite{CTblLib} of ordinary and Brauer character tables.
In particular, \cite{CTblLib} encompasses the data given in the
{\sf Atlas} \cite{Atlas} and the {\sf ModularAtlas} \cite{ModAtlas},
as well as the additional data collected on the 
{\sf ModularAtlasHomepage} \cite{ModAtlasH}.
Data concerning explicit permutation representations,
ordinary and modular matrix representations, and the embedding of
(maximal) subgroups of sporadic simple groups is available in the
{\sf AtlasOfGroupRepresentations} \cite{AtlasGrpRep}, and through the
\GAP{} package {\sf AtlasRep} \cite{AtlasRep}.
To compute with matrix representations over finite fields we use the 
\MA{} \cite{Parker,MA} and its extensions described in \cite{LMR,LuxSz}.

\AbsT{Enumerating long orbits.}
As our computational workhorse to facilitate computations with 
(large) permutation representations we use the \GAP{} package 
\ORB{} \cite{ORB}, whose orbit enumeration techniques are comprehensively
described in \cite{Habil,MNW}. 
For convenience, we give a brief sketch of the approach:

\abs
Let $G$ be a (large) finite group, and let $\cO$ be a (large) 
transitive $G$-set, which we assume to be implicitly given,
for example as a $G$-orbit of a vector $v_1$ in an $F[G]$-module $V$ over 
a finite field $F$. Letting $H\leq G$ be a (still large) subgroup,
we are interested in finding the $H$-orbits $\cO_j\sseq\cO$, their length 
$n_j$, representatives $v_j\in\cO_j$, elements $g_j\in G$ such that
 $v_1\cdot g_j=v_j$, and the point stabilizers $H_j=\Stab_H(v_j)$. 
To achieve this, we assume to be able to compute efficiently within $H$ 
(but not within $G$), for example by having a (smallish) faithful permutation
representation of $H$ at hand.

\abs
To find the $H$-orbits $\cO_j$, we choose a (smallish) helper subgroup 
$K\leq H$, and enumerate the various $\cO_j$ by the
$K$-orbits they contain. To do so, we choose a (not too small)
helper $K$-set $\cQ$ together with a homomorphism $\pi_K\cn\cO\ra\cQ$ 
of $K$-sets, which again we assume to be implicitly given, for example
by an $F[K]$-quotient module of $V$.
We assume that $K$ has sufficiently long orbits in $\cQ$, 
and that we are able to classify them, by giving representatives,
their point stabilizers in $K$, as well as complete Schreier trees.
Thus for the $K$-action on $\cQ$ we are facing a similar problem 
as for the $H$-action on $\cO$, apart from the requirement 
on Schreier trees; so we can just recurse. 

\abs
For any $K$-orbit in $\cQ$, we choose a representative, called its 
`distinguished point'. Then, for any $K$-orbit $\cO'\sseq\cO$, the 
$\pi_K$-preimages of the distinguished point of $\pi_K(\cO')\sseq\cQ$ 
are likewise called the distinguished points of $\cO'$.
Hence to enumerate an $H$-orbit $\cO_j$ by enumerating the $K$-orbits 
it contains, we only have to store the associated distinguished points, 
and a Schreier tree telling us how to reach them from the orbit 
representative $v_j$.

\abs 
For any $\cO_j$ we are content with finding only as many 
$K$-orbits contained in it which are needed to cover (more than) half of it;
this is equivalent to knowing $n_j$ and $|H_j|$. Then we have a 
randomized membership test for $\cO_j$, and a deterministic test to 
decide whether the $\cO_j$ found are actually pairwise disjoint.

\section{Condensing induced modules}\label{indcond}

\AbsT{Endomorphisms of permutation modules.}
We recall some facts about the structure of endomorphism algebras
of permutation modules, thereby fixing the notation used in the 
sequel; as a general reference, see \cite[Ch.II.12]{Landrock}.
 
\abs
Let $G$ be a finite group, let $H\leq G$ be a subgroup, let $\cO$
be a transitive $G$-set with associated point stabilizer $H$, that is 
there is $v_1\in\cO$ such that $\Stab_G(v_1)=H$, and let $n:=|\cO|=[G\cn H]$.

\abs
If $R$ is a principal ideal domain, let $R[\cO]$ be the 
associated permutation $R[G]$-lattice. For subgroups $L\leq M\leq G$ let 
$\Fix_{R[\cO]}(M):=\{v\in R[\cO];vg=v\text{ for all }g\in M\}\leq R[\cO]$ 
be the $R$-sublattice of $M$-fixed points, and let 
$$ \Tr_L^M\cn\Fix_{R[\cO]}(L)\ra\Fix_{R[\cO]}(M)\cn 
   x\mt x\cdot\sum_{g\in L\bsl M}xg $$ 
be the associated trace operator, where $g$ runs through
a set of representatives of the cosets of $L$ in $M$.
For $L=\{1\}$ we just write $\Tr^M:=\Tr_{\{1\}}^M$.

\abs 
Let $E_R:=\End_{R[G]}(R[\cO])$ be the $R$-algebra of $R[G]$-endomorphisms 
of $R[\cO]$. Then $E_R$ is $R$-free of rank $r=|H\bsl G/H|$, that is 
the number of double cosets of $H$ in $G$. In other words, we have 
$r=\lr{1_H^G,1_H^G}_G$, where $1_H^G$ is the permutation character afforded
by $\cO$, and $\lr{\cdot,\cdot}_G$ denotes the usual scalar product on the 
complex class functions on $G$. More precisely:

\abs
Let $\{v_1\ld v_r\}$ be a set of representatives of the $H$-orbits
$\cO_i:=(v_i)^H\sseq\cO$, where $v_1$ is as specified above, 
let $g_i\in G$ such that $v_1g_i=v_i$, let 
$H_i:=\Stab_H(v_i)=H^{g_i}\cap H$, and  let $n_i:=|\cO_i|=[H\cn H_i]$. 
Then $E_R$ has a distinguished $R$-basis $\{A_1\ld A_r\}$, being called its
Schur basis, where $A_i$ is given by 
$$ A_i\cn v_1\mt\cO_i^+:=\sum_{v\in\cO_i}v=v_i\cdot\Tr_{H_i}^H ,$$ 
and extension to all of $\cO$ by $G$-transitivity.

\AbsT{Restriction to subgroups.}
Keeping the above notation, let $U\leq G$ be a subgroup.
Then we may consider $\cO$ as an intransitive $U$-set:

\abs
Let $\{\om_1\ld\om_s\}$ be a set of representatives of the $U$-orbits
$\Om_j:=(\om_j)^U\sseq\cO$, where $\om_1:=v_1$ and $s=|H\bsl G/U|$,
let $\gm_j\in G$ such that $v_1\gm_j=\om_j$, and let 
$U_j:=\Stab_U(\om_j)=H^{\gm_j}\cap U$; then we have $|\Om_j|=[U\cn U_j]$.

\abs
We get a direct sum decomposition $R[\cO]|_U=\bigoplus_{j=1}^s R[\Om_j]$
into transitive permutation $R[U]$-lattices, and a corresponding decomposition 
$$ \cE_R:=\End_{R[U]}(R[\cO]|_U)
  =\bigoplus_{j=1}^s\bigoplus_{k=1}^s \Hom_{R[U]}(R[\Om_j],R[\Om_k]) .$$
We abbreviate $\cE_{jk,R}:=\Hom_{R[U]}(R[\Om_j],R[\Om_k])$. Then 
$\cE_{jk,R}$ has a distinguished $R$-basis $\{\cA_{jk,1}\ld\cA_{jk,t}\}$,
again called its Schur basis, given as follows:

\abs
Let $\{u_{jk,1}\ld u_{jk,t}\}\sseq U$ be a set of representatives of
the double cosets $U_k\bsl U/U_j$, where $t:=|U_k\bsl U/U_j|$, and let
$\om_{jkl}:=\om_k\cdot u_{jkl}\in\Om_k$, for $l\in\{1\ld t\}$.
Then the $U_j$-orbits in $\Om_k$ are given as 
$\Om_{jkl}:=(\om_{jkl})^{U_j}$, where 
$$ \Stab_U(\om_{jkl})=H^{\gm_k u_{jkl}}\cap U
   = H^{\gm_k u_{jkl}}\cap U^{u_{jkl}}
   = (H^{\gm_k}\cap U)^{u_{jkl}} 
   = U_k^{u_{jkl}} ,$$
implying that $|\Om_{jkl}|=\frac{|U_j|}{|U_k^{u_{jkl}}\cap U_j|}$.
Then $\cA_{jkl}\in\cE_{jk,R}$ is given by 
$$ \cA_{jkl}\cn\om_j\mt\Om_{jkl}^+ 
   =\om_{jkl}\cdot\Tr_{U_k^{u_{jkl}}\cap U_j}^{U_j}
   =\om_k\cdot u_{jkl}\cdot\Tr_{U_k^{u_{jkl}}\cap U_j}^{U_j} ,$$
and extension to all of $\Om_j$ by $U_j$-transitivity.  
The action of $u_{jkl}\cdot\Tr_{U_k^{u_{jkl}}\cap U_j}^{U_j}$ 
only depends on the parameters $j,k,l$, but not on the particular choice 
of the $u_{jkl}$.

\AbsT{Embedding endomorphisms.}\label{embed}
Next, we describe the embedding of $E_R$ into $\cE_R$, 
in terms of their Schur bases: 

\abs
For $i\in\{1\ld r\}$ and $j\in\{1\ld s\}$ we have
$$ \om_j\cdot A_i=v_1\gm_j\cdot A_i=v_1\cdot A_i\cdot\gm_j 
   =v_i\cdot\Tr_{H_i}^H\cdot\gm_j=\cO_i^+\cdot\gm_j .$$
Hence, for $k\in\{1\ld s\}$ and $l\in\{1\ld t\}$, where $t=|U_k\bsl U/U_j|$,
let the {\bf $U$-orbit counting numbers} be defined as
$$ c_{jkl}(i):=|\cO_i\gm_j\cap\Om_{jkl}| \quad\text{and}\quad 
c_{jk}(i):=|\cO_i\gm_j\cap\Om_k|=\sum_{l=1}^t c_{jkl}(i) .$$
Recall that $U_j=H^{\gm_j}\cap U$, so that $\cO_i\gm_j$ is $U_j$-stable,
thus so is $\cO_i\gm_j\cap\Om_{jkl}$, implying that $\Om_{jkl}$ is 
either disjoint from $\cO_i\gm_j$, or contained in it, so that either 
$c_{jkl}(i)=0$ or $c_{jkl}(i)=|\Om_{jkl}|$. Thus we have 
$$ A_i=\sum_{j=1}^s\sum_{k=1}^s\sum_{l=1}^{|U_k\bsl U/U_j|} 
       \frac{c_{jkl}(i)}{|\Om_{jkl}|}\cdot\cA_{jkl}\in\cE_R ,$$
where the coefficients are in $\{0,1\}$, saying that $A_i$ splits 
into a sum of certain pairwise distinct Schur basis elements of $\cE_R$.
This is illustrated by the following generic examples:

\abs\bfi
Let $U=G$; thus $\cE_R=E_R$. We have $s=1$ and $\Om_1=\cO$,
so that $\Om_{1,1,l}=\cO_l$ for $l\in\{1\ld r\}$, where 
$|U_1\bsl U/U_1|=|H\bsl G/H|=r$. This yields 
$c_{1,1,l}(i)=\dt_{i,l}\cdot|\cO_l|$ for $i\in\{1\ld r\}$,
and the above triple sum boils down to the tautology $A_i=A_i$.

\abs\bfii
Let $U=\{1\}$; thus $\cE_R=\End_R(R[\cO])\cong R^{n\tm n}$. 
We have $s=n$, and $\Om_k=\{v_1\gm_k\}$ is a singleton set, and 
$|U_k\bsl U/U_j|=1$ for $j,k\in\{1\ld n\}$.
This yields $c_{j,k,1}(i)=1$ if $v_1\gm_k\in\cO_i\gm_j$,
and $c_{j,k,1}(i)=0$ otherwise, for $i\in\{1\ld r\}$.
Hence, identifying $\cA_{j,k,1}$ with $E_{jk}\in R^{n\tm n}$, 
having entry $1$ at position $[j,k]$, and zero entries otherwise,
we recover the natural representation of $E_R$ on $R[\cO]$.

\abs\bfiii
Let $U=H$. We have $s=r$, and $\Om_j=\cO_j$ for $j\in\{1\ld r\}$. Thus 
we get the ($H$-)orbit counting numbers $c_{jk}(i)=|\cO_i g_j\cap\cO_k|$.
These are related to the regular representation of $E_R$ with respect 
to the Schur basis as follows: 

\abs
For $i,j\in\{1\ld r\}$ we write $A_jA_i=\sum_{k=1}^r p_{jk}(i)\cdot A_k$, 
the associated structure constants $p_{jk}(i)\in\N_0$ being called 
intersection numbers. Then we have
$$ v_1\cdot A_jA_i=v_1g_j\cdot\Tr_{H_j}^H\cdot A_i
=v_1\cdot A_i\cdot g_j\cdot\Tr_{H_j}^H =\cO_i^+\cdot g_j\cdot\Tr_{H_j}^H .$$
From $v_1\cdot A_k=\cO_k^+$, and
$|\cO_ig_jh\cap\cO_k| 
=|\cO_ig_j\cap\cO_k|=c_{jk}(i)$,
for $h\in H$, we get
$$ p_{jk}(i)=\sum_{h\in H_j\bsl H}\frac{|\cO_ig_jh\cap\cO_k|}{|\cO_k|} 
 = \frac{n_j}{n_k}\cdot c_{jk}(i) .$$

\AbsT{Condensing $\cE$.}\label{condI}
\bfa
Let $V$ be an $R[U]$-lattice. 
Then $\cE_R$ acts naturally (from the left) on 
$\ccH(V):=\Hom_{R[U]}(R[\cO]|_U,V)$, by 
$$ \al\cn\ccH(V)\ra\ccH(V)\cn\ph\mt\al\cdot\ph,\quad\text{for }\al\in\cE_R .$$
Note that we have
$\ccH(V)\cong\Hom_{R[G]}(R[\cO],V^G)=\Hom_{R[G]}(R_H^G,V^G)$,
where $R_H$ denotes the trivial $R[H]$-module,
and superscripts denote induction. 

\abs
The direct sum decomposition $R[\cO]|_U=\bigoplus_{j=1}^s R[\Om_j]$
entails a decomposition 
$$ \ccH(V)=\bigoplus_{j=1}^s \Hom_{R[U]}(R[\Om_j],V) 
   \cong\bigoplus_{j=1}^s\Fix_V(U_j) ;$$
we write $\ph=\sum_{j=1}^s\ph_j$. 
The latter isomorphism of $R$-lattices is given component-wise by 
$\Hom_{R[U]}(R[\Om_j],V)\ra\Fix_V(U_j)\cn\ph_j\mt\om_j\ph_j$.

\abs
For the Schur basis element $\cA_{jkl}\in\cE_{jk,R}$,
letting $\ph_i\in\Hom_{R[U]}(R[\Om_i],V)$, we have $\cA_{jkl}\cdot\ph_i=0$ 
whenever $i\neq k$. If $i=k$, then we get
$$ \om_j\cdot\cA_{jkl}\ph_k
  =\om_k u_{jkl}\cdot\Tr_{U_k^{u_{jkl}}\cap U_j}^{U_j}\cdot\ph_k
  =\om_k\ph_k\cdot u_{jkl}\cdot\Tr_{U_k^{u_{jkl}}\cap U_j}^{U_j}
\in\Fix_V(U_j), $$
where indeed $\om_k\ph_k\in\Fix_V(U_k)$, hence
$\om_k\ph_k\cdot u_{jkl}\in\Fix_V(U_k^{u_{jkl}})$.
Thus, in terms of fixed spaces, $\cA_{jkl}$ annihilates $\Fix_V(U_i)$
for $i\neq k$, and for $i=k$ we get
$$ \cA_{jkl}\cn\Fix_V(U_k)\ra\Fix_V(U_j)\cn 
   v\mt v\cdot u_{jkl}\cdot\Tr_{U_k^{u_{jkl}}\cap U_j}^{U_j} .$$

\abs\bfb
In view of the application envisaged here, let $R[\cO]$ be a 
projective $R[G]$-module, which is equivalent to $|H|$ being a unit in $R$.
Then let 
$$ e_H:=\frac{1}{|H|}\cdot\Tr^H\in R[H] $$ 
be the associated `fixed-point' idempotent, that is the primitive idempotent
of $R[H]$ associated with the trivial representation of $H$; recall that $e_H$
projects any $R[H]$-lattice onto its $R$-sublattice of $H$-fixed points.
Moreover, we have $R[\cO]\cong R_H^G\cong e_H R[G]$,
so that as $R$-lattices we get 
$$ \ccH(V) = \Hom_{R[U]}(R[\cO]|_U,V)
   \cong \Hom_{R[G]}(e_H R[G],V^G)\cong V^G\cdot e_H .$$

\abs
This shows that $\ccH(V)$ can be seen as the `condensed module' of the 
induced module $V^G$, with respect to the `condensation subgroup' $H$.
A technique to compute condensed modules of shape $V^G\cdot e_H$ for 
subgroups $U$ of smallish index in $G$ has been invented in \cite{IndCond};
the present approach now allows for both subgroups $U$ and $H$ to have
large index. In the spirit of `condensation techniques', recalling that 
$|\Om_{jkl}|=\frac{|U_j|}{|U_k^{u_{jkl}}\cap U_j|}$,
on $\Fix_V(U_k^{u_{jkl}})$ we get
$$ \Tr_{U_k^{u_{jkl}}\cap U_j}^{U_j}
  = e_{U_k^{u_{jkl}}\cap U_j}\cdot\Tr_{U_k^{u_{jkl}}\cap U_j}^{U_j}
  = \frac{|\Om_{jkl}|}{|U_j|}\cdot\Tr^{U_j}
  = |\Om_{jkl}|\cdot e_{U_j} .$$
Thus, in terms of fixed-point idempotents,
the action of $\cA_{jkl}$ can be written as
$$ \cA_{jkl}\cn\Fix_V(U_k)\ra\Fix_V(U_j)\cn 
   v \mt |\Om_{jkl}|\cdot v\cdot e_{U_k}u_{jkl}e_{U_j} .$$ 

\abs
Note that the latter `condensation formula' actually holds more generally 
if $R[\cO]|_U$ is a projective $R[U]$-module, that is all $R[\Om_j]$ are 
projective $R[U]$-modules, which is equivalent to $\prod_{j=1}^s|U_j|$
being a unit in $R$.

\AbsT{Condensing $E$.}\label{condII}
Finally, combining the above observations, 
still assuming that $R[\cO]$ is a projective $R[G]$-module, we derive
a `condensation formula' for the action of the Schur basis elements
of $E_R$ on $\ccH(V)\cong\bigoplus_{j=1}^s\Fix_V(U_j)$:

\abs
Fixing $R$-bases for the fixed spaces $\Fix_V(U_j)$,
the action of $\cE_R$ is given by block matrices, where the blocks
in position $[j,k]$ have size $\rk_R(\Fix_V(U_j))\tm\rk_R(\Fix_V(U_k))$,
for $j,k\in\{1\ld s\}$. Then the matrix representing the Schur basis element
$\cA_{jkl}\in\cE_R$, where $l\in\{1\ld t\}$ and $t:=|U_k\bsl U/U_j|$,
has its only non-zero block in position $[j,k]$, where the latter block 
represents the $R$-linear map $\Fix_V(U_k)\ra\Fix_V(U_j)$ induced by the 
action of $|\Om_{jkl}|\cdot e_{U_k}u_{jkl}e_{U_j}$.

\abs
Let $A_i\in E_R$ be a Schur basis element, where $i\in\{1\ld r\}$.
Then $A_i$ is represented by a block matrix as above, whose block in 
position $[j,k]$ represents the $R$-linear map $\Fix_V(U_k)\ra\Fix_V(U_j)$
induced by the action of
$$ \sum_{l=1}^{|U_k\bsl U/U_j|} c_{jkl}(i)\cdot e_{U_k}u_{jkl}e_{U_j}
  =\sum_{l=1}^{|U_k\bsl U/U_j|} 
   |\cO_i\gm_j\cap\Om_{jkl}|\cdot e_{U_k}u_{jkl}e_{U_j} .$$ 
Since the elements $u_{jkl}\in U$ may be chosen arbitrarily as
representatives of the double cosets $U_k\bsl U/U_j$, 
for any $v\in\Om_k$ we let $u_k(v)\in U$ be any element such that 
$\om_k\cdot u_k(v)=v$. Then the block in position $[j,k]$ represents the map
$$ \sum_{v\in\cO_i\cap(\Om_k\gm_j^{-1})} 
e_{U_k}\cdot u_k(v\gm_j)\cdot e_{U_j} .$$

\section{Characters of $J_4$ and its subgroups}\label{chars}

\Abs
From now on let $G:=J_4$ and $p:=11$. 
Let $\Irr(G)$ be the set of irreducible 
(ordinary) characters of $G$. We order the conjugacy classes of $G$ 
and the $\Irr(G)$ is specified in \cite{Atlas}, thus we may identify 
$\Irr(G)$ with the character table of $G$. 

\abs
The principal $p$-block $B_0$ is the only one 
of positive defect. There are $k_0:=49$ irreducible characters
and $l_0:=40$ irreducible Brauer characters belonging to $B_0$.
According to \cite{ModAtlasH}, this is essentially all what until 
recently has been known about the decomposition numbers of $B_0$.

\abs
Now, in \cite{part1} we have been able to determine four of the 
projective indecomposable characters of $B_0$, amongst them the one 
belonging to the projective cover of the trivial module. The decomposition
of newly found projective indecomposable characters $\Psi_\al$, where 
$\al\in\{1\ld 4\}$, into irreducible characters is reproduced in 
Table \ref{pchtbl}, where $a\in\{0,1\}$.
We also indicate the ordinal numbers of the irreducible characters
occurring, their degree and their character field, where $r_n:=\sqrt{n}$ 
denotes the positive square root of $n\in\N$.

\abs
Actually, $\chi_{19/20}$, $\chi_{23/24}$, $\chi_{36/37}$, and $\chi_{38/39}$ 
form four pairs of mutually algebraically conjugate characters. Since the 
quadratic fields $\Q(r_3)$, $\Q(r_5)$, and $\Q(r_{33})$ are disjoint, there 
are Galois automorphisms of $\Q(r_3,r_5,r_{11})$ inducing each of the 
involutions $(\chi_{19}\Lra\chi_{20})$, $(\chi_{23}\Lra\chi_{24})$, and 
$(\chi_{36}\Lra\chi_{37})(\chi_{38}\Lra\chi_{39})$. Thus we may choose the 
rows of the decomposition matrix belonging to $\chi_{23/24}$ 
and $\chi_{36/37}$ as is shown in Table \ref{pchtbl}, while the 
non-$p$-rational characters $\chi_{19/20}$ have the same restriction 
to the $p$-regular conjugacy classes of $G$ anyway.
But then, as we will see below, there is no further choice possible for 
$\chi_{38/39}$, leaving two possible cases parameterized by $a\in\{0,1\}$. 

\begin{table}\caption{Projective indecomposable characters of $G$,
                      taken from \cite{part1}.}\label{pchtbl}
$$ \begin{array}{|r|r|r||rrrr|} \hline
\chi & \chi(1) & \Q(\chi) & \Psi_1 & \Psi_2 & \Psi_{3,a} & \Psi_{4,a} \\
\hline \hline
 1 &          1 &            &  1 & . & . & . \\
 8 &     889111 &            &  . & 1 & . & . \\
11 &    1776888 &            &  . & . & 1 & . \\
14 &    4290927 &            &  1 & . & . & . \\
19 &   35411145 & \Q(r_{33}) &  1 & . & 1 & . \\
20 &   35411145 & \Q(r_{33}) &  1 & . & 1 & . \\
21 &   95288172 &            &  1 & . & 1 & . \\
22 &  230279749 &            &  1 & . & . & . \\
23 &  259775040 &    \Q(r_3) &  . & . & . & 1 \\
24 &  259775040 &    \Q(r_3) &  . & 1 & . & . \\
29 &  460559498 &            &  . & . & . & 1 \\
30 &  493456605 &            &  . & . & . & 1 \\
32 &  786127419 &            &  . & 1 & . & . \\
36 &  885257856 &    \Q(r_5) &  1 & . & . & . \\
37 &  885257856 &    \Q(r_5) &  . & . & . & 1 \\
38 & 1016407168 &    \Q(r_5) &  . & . & a &1-a\\
39 & 1016407168 &    \Q(r_5) &  . & . &1-a& a \\
51 & 1842237992 &            &  . & . & . & 1 \\
\hline \end{array} $$
\absr\end{table}

\Abs
In order to get a comprehensive overview about the possible choices
on the character theoretic side, and what has to be decided explicitly
in the end, we use the following terminology:
Let $\cA(G):=\Aut(\Irr(G))$ be the group
of table automorphisms of $\Irr(G)$, that is the permutations of the 
conjugacy classes of $G$ compatible with power maps and inducing 
permutations of the rows of $\Irr(G)$.

\abs
According to \cite{CTblLib}, the group $\cA(G)$ has order $432$, 
is generated by 
$$ \begin{array}{rcl} 
\cA(G)&\hspace*{-0.5em}=&\hspace*{-0.5em}\lr{
(12,13)(24,25)(26,27)(32,33)(39,40)(48,49)(55,56),\, (43,44,45),\\ 
&&(30,31)(53,54),\, (37,38),\; (46,47)(61,62),\, (50,51,52),\, (57,58,59)} .\\
\end{array} $$
Its action on $\Irr(G)$ is given as
$$ \begin{array}{rcl} 
\cA(G)&\hspace*{-0.5em}\ra&\hspace*{-0.5em}\lr{
 (2,3)(4,5)(6,7)(9,10)(12,13)(15,16)(17,18),\, (19,20)(33,34), \\
&& (23,24),\, (36,37)(38,39),\, (46,47,48),\, (53,54,55),\, (56,57,58)} .\\
\end{array} $$
In particular, the latter contains the action of the Galois automorphisms 
mentioned above, where on the irreducible characters considered here we 
indeed see the involution $(19,20)$. (On all of $\Irr(B_0)$ we see
$(19,20)(33,34)$ instead, which has to be taken into account as soon as as 
projective characters also having constituents $\chi_{33/34}$ are considered.)

\Abs
In the sequel, for various subgroups $M<G$, we will compute the set of
`possible class fusions' from the conjugacy classes of $M$ to those of $G$, 
that is the maps compatible with power maps and restrictions of
irreducible characters. This set is acted on naturally by 
$\cA(M)\tm\cA(G)$ via $[\al,\bt]\cn f\mt\al^{-1}\cdot f\cdot\bt$.
Possible class fusions are considered equivalent if they belong to the same 
$(\cA(M)\tm\cA(G))$-orbit. Hence `choosing' a class fusion amounts to 
picking an orbit representative, $f$ say, and keeping it fixed. But 
this restricts the table automorphisms remaining admissible in subsequent 
`choices up to equivalence' to $\Stab_{\cA(M)\tm\cA(G)}(f)\leq\cA(M)\tm\cA(G)$,
and its projections $\cA_f(M)\leq\cA(M)$ and $\cA_f(G)\leq \cA(G)$ 
to the first and second direct factors, respectively.

\abs
Now we bring Brauer characters into play: Let $\IBr_p(G)$ be the 
(as yet unknown) $p$-modular Brauer character table of $G$, whose 
columns are identified with the $p$-regular conjugacy classes of $G$.
This amounts to saying that the class fusion from $p$-regular 
to all conjugacy classes of $G$ has been chosen, so that we have 
to go down to the admissible subgroup $\cA_p(G)\leq\cA(G)$ inducing 
permutations of the rows of $\IBr_p(G)$.

\abs
Similarly, let $\IPr_p(G)$ be the (as yet unknown) table of $p$-modular 
projective indecomposable characters of $G$. Since $\IPr_p(G)$ is the
dual basis of $\IBr_p(G)$, extended by zeroes on the $p$-singular conjugacy
classes of $G$, with respect to the usual scalar product on the complex class 
functions on $G$, the group $\cA_p(G)$ coincides with the subgroup 
of $\cA(G)$ inducing permutations of the rows of $\IPr_p(G)$.
If only an $\cA_p(G)$-stable subset $\Psi\sseq\IPr_p(G)$ is known, 
then the subgroup $\cA_\Psi(G)\leq\cA(G)$ inducing permutations of $\Psi$ 
contains $\cA_p(G)$, thus can serve as an upper approximation of the latter.

\Abs
Now, let $H<G$ be a (maximal) subgroup of shape $2^{10}\cn L_5(2)$.
Then $H$ is an $11'$-subgroup of order $|H|=10\,239\,344\,640$.
It turns out that $\cA(H)$ has order $24$, and that there are six
possible class fusions from $H$ to $G$, consisting of a single 
$(\cA(H)\tm\cA(G))$-orbit.
As a representative $f$ we choose the class fusion stored in \cite{CTblLib}.
We get $\cA_f(G)=\cA(G)$, saying that upon choosing the class fusion 
from $H$ to $G$ all table automorphisms of $G$ remain admissible.

\abs
Moreover, let $1_H^G$ be the permutation character afforded by a
transitive $G$-set with associated point stabilizer $H$. Then $1_H^G$ 
is $\cA(G)$-invariant, and it is a projective character of $G$, which 
by \cite{part1} splits into four projective indecomposable characters as  
$1_H^G=\Psi_1+\Psi_2+\Psi_{3,a}+\Psi_{4,a}$. 
Hence the set $\Psi:=\{\Psi_1,\Psi_2,\Psi_{3,a},\Psi_{4,a}\}$ 
is $\cA_p(G)$-stable.
From the description of the action of $\cA(G)$ on $\Irr(G)$ we infer
that $\cA_{\Psi}(G)$ fixes $\Psi$ element-wise, and thus so does $\cA_p(G)$. 
Hence we get $\cA_\Psi(G)=\Stab_{\cA(G)}(\chi_{23},\chi_{36})$,
for both $a\in\{0,1\}$, that is
$$ \begin{array}{rcl} 
\cA_\Psi(G)&\hspace*{-0.5em}=&\hspace*{-0.5em}\lr{
(12,13)(24,25)(26,27)(32,33)(39,40)(48,49)(55,56), \\
&& (43,44,45),\, (46,47)(61,62),\, (50,51,52),\, (57,58,59)} .\\
\end{array} $$
Thus we have $[\cA(G)\cn\cA_\Psi(G)]=4$, reflecting the couple of choices
made between two alternatives each, and $\cA_\Psi(G)$ acts on $\Irr(G)$ as 
$$ \begin{array}{rcl} 
\cA_\Psi(G)&\hspace*{-0.5em}\ra&\hspace*{-0.5em}\lr{
 (2,3)(4,5)(6,7)(9,10)(12,13)(15,16)(17,18), \\
&& (19,20)(33,34),\, (46,47,48),\, (53,54,55),\, (56,57,58)} .\\
\end{array} $$
We conclude that $\cA_\Psi(G)$ acts on the constituents of 
the permutation character $1_H^G$ as a subgroup of
$\lr{(19,20)}$, and so does $\cA_p(G)$. In particular, the cases
$a\in\{0,1\}$ are genuinely different, so that we have to decide  
which one holds.

\Abs
In order to do so, let $U<G$ be a (maximal) subgroup of shape $U_3(11)\cn 2$,
having order $|U|=141\,831\,360$, 
and let $U'\cong U_3(11)$ be its derived subgroup of index $2$.
Moreover, let $S<U'$ be a Sylow $11$-subgroup, hence $S$ is extra-special 
of shape $11^{1+2}_+$, and is a Sylow $11$-subgroup of $G$ as well.
 
\abs
Let $\IBr_p(U)$ be as specified in \cite{ModAtlas}.
In particular, let $S_8$ be the (absolutely irreducible) adjoint module
of $U'$ of degree $8$, and let $S_8^{\pm}$ be its extensions to $U$, 
where the Brauer character of $S_8^{\pm}$ on the conjugacy class of 
involutions not contained in $U'$ has value $\pm 2$. Then $S_8^+$ and
$S_8^-$ have Brauer characters $\ph_3$ and $\ph_4$, respectively.
Moreover, let $\Phi_8^{\pm}$ be the projective indecomposable characters
of $U$ associated with $S_8^{\pm}$.

\abs
It turns out that the group $\cA(U)$ has order $96$,
but the admissible subgroup $\cA_p(U)\leq\cA(U)$ has order $2$, 
whose non-trivial element is the transposition interchanging 
the conjugacy classes of elements of order $44$ not belonging to $U'$. 
It turns out that there are $24$ possible class fusions from $U$ to $G$, 
which fall into two $(\cA_p(U)\tm\cA(G))$-orbits of length $12$;
the latter are even $\cA(U)$-invariant. Orbit representatives are given 
(in terms of conjugacy class numbers) as
$$ \begin{array}{c}  
 [ 1, 2, 4, 5, 6, 8, 8, 10, 14, 17, 17, 19, 20, 20, 21, 22, 22, \\
   30, 31, 34, 50, 51, 52, 50, 51, 52, 53, 54, 53, 54, 60, \\
   3, 5, 11, 15, 18, 18, 21, 30+y, 31-y, 35, 37, 38, 60, 60 ],
\end{array} $$
where $y\in\{0,1\}$, the case $y=0$ being the one stored in \cite{CTblLib}. 
The conjugacy classes of $U$ whose fusion to $G$ depends on $y$ 
are those containing elements of order $20$ not belonging to $U'$.

\abs
The picture changes when we restrict to $\cA_\Psi(G)$:
Both of the above orbits split into four $(\cA_p(U)\tm\cA_\Psi(G))$-orbits 
of length three. Thus now there are eight orbits, 
representatives of which are given by the maps
$$ \begin{array}{c}  
 [ 1, 2, 4, 5, 6, 8, 8, 10, 14, 17, 17, 19, 20, 20, 21, 22, 22, \\ 
   30+x, 31-x, 34, 50, 51, 52, 50, 51, 52, 53+x, 54-x, 53+x, 54-x, 60, \\
   3, 5, 11, 15, 18, 18, 21, 30+y, 31-y, 35, 37+z, 38-z, 60, 60 ],
\end{array} $$
where $x,y,z\in\{0,1\}$; 
the case $x=y=z=0$ being the one stored in \cite{CTblLib}.
The conjugacy classes of $U$ whose fusion to $G$ depends on the parameters
$x$ or $y$ consist of elements of order $20$ and $40$, where the conjugacy
classes $\{53,54\}$ of $G$ square to the conjugacy classes $\{30,31\}$;
the conjugacy classes of $U$ whose fusion to $G$ depends on the parameter $z$ 
consist of elements of order $24$.

\Abs\label{max7}
We can do slightly better, as far as the fusion from conjugacy classes
of $U$ consisting of elements of order $20$ to $G$ is concerned:

\abs
To this end, let $N:=N_G(S)<G$, which is a (maximal) subgroup of shape 
$S\cn T\cong 11^{1+2}_+\cn(5\tm 2.\cS_4)$, having order $|N|=319\,440$,
where $T\cong 5\tm 2.\cS_4$ is unique up to $N$-conjugacy. 
Hence we may assume that $T\cap U\cong 5\tm\text{QD}_{16}$, thus 
$$ N_U(S)=N\cap U=S\cn(T\cap U)\cong 11^{1+2}_+\cn(5\tm\text{QD}_{16}) .$$

\abs
We compare the embeddings $T\cap U<U$ and $T\cap U<T$:
It turns out that there are eight possible class fusions from 
$T\cap U$ to $T$, and two possible class fusions from $T$ to $G$. 
Then composition yields two possible class fusions from $T\cap U$ to $G$ 
which factor through $T$.
Similarly, there are four possible class fusions from $T\cap U$ to $U$. 
Then composing either of the $24$ possible class fusions from 
$U$ to $G$ with the latter, and checking whether a possible class 
fusion factoring through $T$ is obtained, leaves $12$ possible class 
fusions from $U$ to $G$. 

\abs
It turns out that these fall into four $(\cA_p(U)\tm\cA_\Psi(G))$-orbits, 
which are given by the parameters $x=y$, leaving the following four
maps from the above list:
$$ \begin{array}{c}  
 [ 1, 2, 4, 5, 6, 8, 8, 10, 14, 17, 17, 19, 20, 20, 21, 22, 22, \\ 
   30+y, 31-y, 34, 50, 51, 52, 50, 51, 52, 53+y, 54-y, 53+y, 54-y, 60, \\
   3, 5, 11, 15, 18, 18, 21, 30+y, 31-y, 35, 37+z, 38-z, 60, 60 ].
\end{array} $$


\Abs\label{mult} 
Having this in place, restricting $\Psi_\al$ to $U$,
for $\al\in\{1\ld 4\}$ and both cases $a\in\{0,1\}$, 
and using the various class fusions for $y,z\in\{0,1\}$, 
we may write $\Psi_\al|_U$ uniquely as an integral linear combination of
projective indecomposable characters of $U$. By construction, the 
multiplicities occurring are $(\cA_p(U)\tm\cA_\Psi(G))$-invariant. 
It turns out that in all cases these multiplicities are non-negative,
so that this does not yield further immediate restrictions. 

\abs
But, amongst others, the multiplicity $[\Psi_\al|_U\cn\Phi_8^\pm]$ 
of $\Phi_8^\pm$ in a direct sum decomposition of $\Psi_\al|_U$ subtly 
depends on the parameters $a,y,z$. We get the following pattern,
where we also indicate $[1_H^G|_U\cn\Phi_8^\pm]$:
$$ \begin{array}{|l||ll|} \hline &\Phi_8^+ & \Phi_8^- \\ 
\hline \hline 
\Psi_1|_U     & 77-y    & 71+y \\
\Psi_2|_U     & 67+z    & 52-z \\
\Psi_{3,a}|_U & 80      & 56+(-1)^{y+a} \\
\Psi_{4,a}|_U & 299+y-z & 260-y+z-(-1)^{y+a} \\
\hline  \hline 
1_H^G|_U      & 523     & 439 \\
\hline 
\end{array} $$


\abs
In order to determine the multiplicities $[\Psi_\al|_U\cn\Phi_8^\pm]$
explicitly, let $\F=\F_{11}$, and let $P_\al$ be the projective
indecomposable $\F[G]$-module affording $\Psi_\al$, where since 
$[1_H^G\cn\Psi_\al]=1$ implies that $\Psi_\al$ is indeed realizable over $\F$. 
We have 
$$ [\Psi_\al|_U\cn\Phi_8^\pm]=\dim_\F(\Hom_{\F[U]}(P_\al|_U,S_8^\pm)) .$$
Using the notation introduced in Section \ref{indcond}, in particular letting 
$\cO$ be a $G$-set affording the permutation character $1_H^G$, we have 
$\F[\cO]\cong\bigoplus_{i=1}^4 P_\al$ as $\F[G]$-modules, 
where $P_\al\cong\F[\cO]\cdot e_\al$ for a set $\{e_1\ld e_4\}\sseq E_\F$ of 
pairwise orthogonal primitive idempotents. 
The natural action of $E_\F\sseq\cE_\F$ on 
$$ \ccH(S_8^\pm):=\Hom_{\F[U]}(\F[\cO]|_U,S_8^\pm)\cong 
\bigoplus_{j=1}^s\Fix_{S_8^\pm}(U_j) $$ 
yields $\Hom_{\F[U]}(P_\al|_U,S_8^\pm)\cong e_\al\cdot\ccH(S_8^\pm)$.
Thus we have to determine the action of $E_\F$ on $\ccH(S_8^\pm)$,
and the $\F$-dimension 
$d_\al:=\dim_\F(e_\al\cdot\ccH(S_8^\pm))$, for $\al\in\{1\ld 4\}$.
At this stage, we switch to explicit computations:

\section{Enumerating $\cO$ again}\label{enum}

\Abs
We pick the $112$-dimensional absolutely irreducible representation 
of $G$ over $\F_2$ from \cite{AtlasRep}, and let $V\cong\F_2^{112}$ be the
underlying module. The representation is given in terms of (two) standard 
generators, in the sense of \cite{Wilson}. Words in the generators providing 
(non-standard) generators of maximal subgroups $2^{10}\cn L_5(2)\cong H<G$
and $U_3(11)\cn 2\cong U<G$ are available in \cite{AtlasRep} as well.

\abs
It turns out that $H$ possesses a $1$-dimensional fixed space in $V$. 
Hence letting $v_1\in V$ be the unique non-zero $H$-fixed vector, we let 
$\cO:=(v_1)^G\sseq V$, providing an implicit realization of $\cO$. 
In \cite{part1}, using \ORB{}, we have already enumerated $\cO$ by
$H$-orbits, of which there are $r:=\lr{1_H^G,1_H^G}_G=27$. 
For the $H$-orbits $\cO_i\sseq\cO$ we have in particular determined
their lengths $n_i$, which are reproduced in Table \ref{orbtbl}, as well as 
elements $g_i\in G$ yielding orbit representatives $v_i:=v_1g_i\in\cO_i$.

\abs
We are now going to enumerate $\cO$ again, this time by $U$-orbits. 
It turns out that there are $s:=\lr{1_H^G,1_U^G}_G=131$ such $U$-orbits;
note that $s$ is independent of the class fusions chosen.
To do so, we set up a new framework to apply \ORB{}, 
adjusted to our present needs. In particular, comparing with \cite{part1}, 
we have to pick another helper subgroup, since the one chosen there 
is a subgroup of $H$, but is not conjugate to a subgroup of $U$.

\begin{table}\caption{$H$-orbit lengths in $\cO$.}\label{orbtbl}
$$ \begin{array}{|r|r|} \hline 
 i & n_i \\ \hline \hline
 1 & 1 \\
 2 & 31 \\
 3 & 930 \\
 4 & 17360 \\
 5 & 26040 \\
 6 & 27776 \\
 7 & 416640 \\
 8 & 416640 \\
 9 & 624960 \\
\hline \end{array} \quad \quad \quad
\begin{array}{|r|r|} \hline
 i & n_i \\ \hline \hline
10 & 333120 \\
11 & 4999680 \\
12 & 6666240 \\
13 & 6666240 \\
14 & 9999360 \\
15 & 13332480 \\
16 & 53329920 \\
17 & 66060288 \\
18 & 79994880 \\
\hline \end{array} \quad \quad \quad
\begin{array}{|r|r|} \hline
 i & n_i \\ \hline \hline
19 & 79994880 \\
20 & 159989760 \\
21 & 159989760 \\
22 & 319979520 \\
23 & 341311488 \\
24 & 1279918080 \\
25 & 1279918080 \\
26 & 2047868928 \\
27 & 2559836160 \\
\hline \end{array} $$
\absr\end{table}

\Abs
By a random search we replace the non-standard generators of $U$ 
we have so far by standard ones. A faithful permutation representation 
of $U$ on $1332$ points, in terms of standard generators, is available in
\cite{AtlasRep}. Actually, the associated point stabilizers are conjugate
to $11^{1+2}_+\cn(5\tm\text{QD}_{16})\cong N_U(S)<U$, which was already 
encountered in \ref{max7}.

\abs
We choose a (single) helper subgroup $K<U$:
Let $z\in U'$ be an involution, which is unique up to $U'$-conjugacy;
specifically, we choose $z$ as the square of the second standard generator 
of $U$, which has order $4$. Then we let 
$$ K:=C_U(z)\cong (\SL_2(11)\tm V_4)\cn 2 ,$$
being computed in the permutation representation of $U$;
we have $|K|=10560$. 

\abs
Keeping $K$ fixed, we choose (two) helper $K$-sets:
It turns out that the restriction of $V$ to $K$ decomposes as 
$V|_K \cong V_{80}\oplus V_{32}$, subscripts denoting $\F_2$-dimension,
where $K$ acts faithfully on both summands. Moreover, $V_{32}$ 
has a unique $\F_2[K]$-quotient $V_{20}$ of $\F_2$-dimension $20$,
on which $K$ acts non-faithfully by its quotient $(L_2(11)\tm V_4)\cn 2$.
Actually, $V_{20}$ is uniserial with a unique (absolutely irreducible) 
constituent of $\F_2$-dimension $10$, on which $K$ acts as $L_2(11)\cn 2$.
As helper $K$-sets we now choose the natural epimorphisms 
$V|_K \ra V_{32} \ra V_{20}$ of $\F[K]$-modules.

%
%

\Abs
We are now prepared to run \ORB{}, in order to find the decomposition 
$\cO=\coprod_{j=1}^{131}\Om_j$ into $U$-orbits:
We randomly choose elements $g\in G$, and check whether $v_1g\in\cO$ 
belongs to one of the $U$-orbits already found. If not, then we have found
a previously unseen $U$-orbit, $\Om_j$ say. In this case, we store 
$\gm_j:=g$ and the orbit representative $\om_j:=v_1\gm_j\in\Om_j$, 
we enumerate half of $\Om_j$, and using the faithful permutation 
representation of $U$ we determine $U_j:=\Stab_U(\om_j)$. 
The lengths of the $U$-orbits in $\cO$ are summarized in 
Table \ref{Uorbitlengthstbl}. Recalling that $|U|=141\,831\,360$,
we infer that the point stabilizers $U_j$ have order at most $360$.

\abs
To detect all $U$-orbits in $\cO$ we need approximately half an hour 
on a single 3 GHz CPU. Setting up \ORB{} anew, and using the $\gm_j$
instead of a random search, the orbit enumeration database is rebuilt 
in about ten minutes of CPU time.
The statistics provided by \ORB{} shows that the `saving factor', 
that is the quotient between the number of points in a $U$-orbit 
actually stored in the orbit enumeration database, and the length
of the piece of the $U$-orbit enumerated, varies between $8022$ 
(for one of the shorter $U$-orbits) and $10497$ (for one of the
regular $U$-orbits). Thus we achieve an average saving factor of 
$\sim 10304 \sim 0,97 \cdot |K|$. The total memory usage of the 
orbit enumeration database amounts to manageable $\sim 115$ MB, 
and the infrastructure needs additional $\sim 195$ MB.

\begin{table}\caption{$U$-orbit lengths in $\cO$.}\label{Uorbitlengthstbl}
$$ \begin{array}{|r|} \hline
393976 \\
738705 \\
984940 \\
2 \tm 1181928 \\
1477410 \\
\hline \end{array} \quad 
\begin{array}{|r|} \hline
2216115 \\
3 \tm 2954820 \\
3939760 \\
5909640 \\
9 \tm 8864460 \\ 
\hline \end{array} \quad 
\begin{array}{|r|} \hline
4 \tm 11819280 \\
14183136 \\
12 \tm 17728920 \\
3 \tm 23638560 \\
2 \tm 28366272 \\
\hline \end{array} \quad 
\begin{array}{|r|} \hline
22 \tm 35457840 \\
2 \tm 47277120 \\
28 \tm 70915680 \\
36 \tm 141831360 \\ \\ 
\hline \end{array} $$
\absr\end{table}

\section{Condensing the induced module $(S_8^\pm)^G$}\label{conclusion}

\AbsT{Finding idempotents.}
We proceed to determine pairwise orthogonal primitive idempotents 
$\{e_1\ld e_4\}\sseq E_\F$, and their action on $\ccH(S_8^\pm)$.
The idea pursued here is inspired by \cite{LMR}.

\abs
In \cite{part1} we have computed the $11$-modular character 
table of $E_\F$; it is reproduced in Table \ref{modtbl}. Here, notation
is chosen such that the irreducible character $\ph_\al$ of $E_\F$ 
corresponds to the projective indecomposable character $\Psi_\al$ of $G$. 
Since all irreducible characters are linear, the character values coincide
with the eigenvalues of the action of the Schur basis elements in the various
irreducible representations. We consider the action of $A_2$, where we 
observe that the character values $[\ph_1(A_2)\ld\ph_4(A_2)]=[9,5,10,1]$ 
are pairwise different. 

\abs
Let $\mu:=\prod_{\al=1}^4(X-\ph_\al(A_2))^{h_\al}\in\F[X]$ 
be the minimum polynomial of the action of $A_2$ in the (faithful)
regular representation of $E_\F$. The multiplicities $h_\al\in\N$ are not 
needed explicitly in the sequel, but they are actually easily determined:
The intersection matrices of $E_\Z$
have been determined in \cite{part1}, so that the minimum polynomial of
their $11$-modular reduction is straightforwardly computed, yielding 
$[h_1\ld h_4]=[5,3,5,4]$.

\abs
Let $\mu_\al:=(X-\ph_\al(A_2))^{h_\al}$, and let 
$\mu'_\al:=\frac{\mu}{\mu_\al}$ be the associated co-factor.
Then, $\mu_\al$ and $\mu'_\al$ being coprime, there are 
$f_\al,f'_\al\in\F[X]$ such that $f_\al\mu_\al+f'_\al\mu'_\al=1\in\F[X]$.
Hence we have $E_\F\cong\bigoplus_{\al=1}^4\ker(\mu_\al(A_2))$ 
as $\F[A_2]$-modules. Moreover, 
$$ e_\al:=f'_\al(A_2)\mu'_\al(A_2)=1-f_\al(A_2)\mu_\al(A_2)\in\F[A_2] $$
annihilates $\ker(\mu'_\al(A_2))=\bigoplus_{\bt\neq\al}\ker(\mu_\bt(A_2))$,
while it acts as the identity on $\ker(\mu_\al(A_2))$.
Thus $\{e_1\ld e_4\}\sseq E_\F$ is a set of pairwise orthogonal, hence 
primitive idempotents. Moreover, $e_\al$ acts on any $E_\F$-module
as a projection onto the generalized eigenspace of $A_2$ with respect to 
the eigenvalue $\ph_\al(A_2)$. In particular, it does so on the 
simple $E_\F$-modules, so that $e_\al$ is associated with the 
irreducible character $\ph_\al$ indeed.

\begin{table}\caption{The character table of $E_\F$.}\label{modtbl}
$$ \begin{array}{|r||rrrrrrrrrrrrrr|} \hline
\ph_\al & 1& 2& 3& 4& 5& 6& 7& 8& 9& 10& 11& 12& 13& 14 \\
\hline
1 & 1&  9& 6& 2&  3& 1& 4& 4&  6& 10& 4& 9& 9& 8 \\
2 & 1&  5& 5& 1&  8& 8& 5& 5& 10&  8& 1& 2& 2& 5 \\
3 & 1& 10& 3& 4& 10& 4& 7& 7&  5&  7& 5& 3& 3& 6 \\
4 & 1&  1& 3& 9&  0& 8& 0& 0&  5&  8& 2& 4& 4& 4 \\
\hline \end{array} $$ \vspace*{0.5em}
$$ \begin{array}{|r||rrrrrrrrrrrrr|} \hline
\ph_\al & 15& 16& 17& 18& 19& 20& 21& 22& 23& 24& 25& 26& 27 \\ 
\hline
1 & 7& 6& 8& 9& 9&  7&  7&  3&  1&  1& 1&  6& 2 \\
2 & 3& 5& 5& 5& 5&  6&  3& 10&  1&  0& 1&  3& 8 \\ 
3 & 5& 1& 1& 7& 7& 10&  9&  0&  4&  3& 6& 10& 5 \\
4 & 5& 9& 0& 3& 3&  2& 10&  8& 10& 10& 1&  0& 0 \\
\hline \end{array} $$
\absr\end{table}

\Abs\label{dim}
Thus we are left with determining the action of $A_2$ on $\ccH(S_8^\pm)$, 
which is given in terms of fixed spaces by the `condensation formula' in 
\ref{condII}. In order to apply it, we first observe that the point 
stabilizers $U_j$ are small enough such that the action of the fixed-point
idempotent $e_{U_j}$ on $S_8^\pm$ is straightforwardly computed by running
through all the elements of $U_j$ explicitly.

\abs
Starting with $v_2\in\cO_2$, the $H$-orbit $\cO_2$,
having length $n_2=31$, is easily enumerated explicitly.
Then, running through the points $v\in\cO_2$,
we apply the elements $\gm_j\in G$ in turn, for $j\in\{1\ld s\}$, and 
check using \ORB{} to which $U$-orbit the point $w:=v\gm_j\in\cO$ belongs. 
If $w\in\Om_k$, say, then we use the functionality \ORB{} readily offers 
to find an element $u_k(w)\in U$ such that $\om_k\cdot u_k(w)=w$.

\abs
Having this in place, we apply the `condensation formula' to compute the 
action of $A_2$ on $\ccH(S_8^\pm)\cong\bigoplus_{j=1}^s\Fix_{S_8^\pm}(U_j)$.
This straightforwardly yields the $\F$-dimension $d_\al^\pm$ of 
the generalized eigenspace of the action of $A_2$ with respect
to the eigenvalue $\ph_\al(A_2)$, and for comparison the multiplicity 
$h_\al^\pm$ of the irreducible factor $X-\ph_\al(A_2)$ in the 
minimum polynomial of this action:
$$ \begin{array}{|r|r||rr|rr|} \hline 
 & h_\al & d_\al^+ & h_\al^+ & d_\al^- & h_\al^- \\ \hline \hline 
\ph_1(A_2) & 5 &  76 & 4 &  72 & 2 \\
\ph_2(A_2) & 3 &  67 & 3 &  52 & 2 \\
\ph_3(A_2) & 5 &  80 & 4 &  55 & 2 \\
\ph_4(A_2) & 4 & 300 & 4 & 260 & 4 \\
\hline \hline
\ccH(S_8^\pm)   &   & 523 &   & 439 & \\
\hline \end{array} $$

\Abs
We just remark that the same approach works for all the shorter 
$H$-orbits in $\cO$, in particular including $\cO_4$ of length $n_4=17360$.

\abs
Now, in \cite{part1} we have shown that $E_\Q=\Q[A_2,A_4]$,
by writing all the Schur basis elements of $E_\Q$ explicitly 
as words in the generating set $\{A_2,A_4\}$.
Letting $\Z_{(11)}\sseq\Q$ be the ring of $11$-adic integers in $\Q$,
it turns out that the latter words actually belong to $\Z_{(11)}[A_2,A_4]$. 
This implies that the Schur basis elements of $E_\F$ are given by 
the very words, now considered via $11$-modular reduction as belonging to 
$\F[A_2,A_4]$. Thus we have $E_\F=\F[A_2,A_4]$.

\abs
As already follows from comparing $h_\al^\pm$ with $h_\al$,
the algebra $E_\F$ acts non-faithfully on $\ccH(S_8^\pm)$,
where it turns out that on $\ccH(S_8^+)$ and $\ccH(S_8^-)$ 
it acts by an algebra of $\F$-dimension $24$ and $15$, 
respectively, while $\dim_\F(E_\F)=27$.

\AbsT{Conclusion.} 
We have found the multiplicities
$[\Psi_\al|_U\cn\Phi^\pm]=d_\al^\pm$ as given in \ref{dim}.
Comparing with the possible parameter choices left in \ref{mult}, we get
$$ y=1 \quad \text{and} \quad z=0 \quad \text{and} \quad a=0 .$$

\abs
Thus the projective indecomposable summands $\Psi_\al$ of the 
permutation character $1_H^G$ are, up to admissible table automorphisms, 
as shown in Table \ref{pchtbl}, upon specifying $a:=0$.
Moreover, the class fusion from $U$ to $G$ is given,
again up to admissible table automorphisms, as follows, 
where it turns out that it differs from the one stored in \cite{CTblLib} 
precisely in the (eight) positions printed in bold face:
$$ \begin{array}{c}  
 [ 1, 2, 4, 5, 6, 8, 8, 10, 14, 17, 17, 19, 20, 20, 21, 22, 22, \\ 
   {\bf 31}, {\bf 30}, 34, 50, 51, 52, 50, 51, 52,
   {\bf 54}, {\bf 53}, {\bf 54}, {\bf 53}, 60, \\
   3, 5, 11, 15, 18, 18, 21, {\bf 31}, {\bf 30}, 35, 37, 38, 60, 60 ]. 
\end{array} $$


\absr


\absr\abs\abs
{\sc
Chair for Algebra and Number Theory, RWTH Aachen University \\
Pontdriesch 14/16, D-52062 Aachen, Germany} \\
{\sf e-mail: juergen.mueller@math.rwth-aachen.de}

\end{document}